\documentclass[12pt,english]{extarticle}
\usepackage[T1]{fontenc}
\usepackage[latin9]{inputenc}
\usepackage{geometry}
\usepackage{babel}
\usepackage{amsmath}
\usepackage{amsthm}
\usepackage{amssymb}
\usepackage[unicode=true,pdfusetitle,
 bookmarks=true,bookmarksnumbered=false,bookmarksopen=false,
 breaklinks=false,pdfborder={0 0 1},backref=false,colorlinks=false]
 {hyperref}

\makeatletter
\numberwithin{equation}{section}
\numberwithin{figure}{section}
\theoremstyle{plain}
\newtheorem{thm}{\protect\theoremname}
  \theoremstyle{remark}

\usepackage{color}

\makeatother

  \providecommand{\remarkname}{Remark}
\providecommand{\theoremname}{Theorem}

\begin{document}

\title{A high-order corrector estimate for a semi-linear elliptic
system in perforated domains}

\author{Vo Anh Khoa\thanks{Author for correspondence. Mathematics and Computer Science Division,
Gran Sasso Science Institute, L'Aquila, Italy. (khoa.vo@gssi.infn.it,
vakhoa.hcmus@gmail.com)}}

\maketitle
\begin{abstract}
We derive in this note a high-order corrector estimate for the homogenization of a
microscopic semi-linear elliptic system posed in perforated domains.
The major challenges are the presence
of nonlinear volume and surface reaction rates. This type of correctors justifies mathematically the convergence rate
of formal asymptotic expansions for the two-scale homogenization
settings. As main tool, we follow the standard approach by the energy-like
method to investigate the error estimate between the micro and macro
concentrations and micro and macro concentration gradients.
This work aims at generalizing the results reported in \cite{CP99,KM16}.
\end{abstract}

\section{Introduction}

This note is devoted to the derivation of a high-order corrector for a microscopic
semi-linear elliptic system posed in heterogeneous/perforated domains.
In the terminology of homogenization, a \emph{corrector} or \emph{corrector
estimate} wants to quantify the error between the approximate
solution (governed by the asymptotic procedure) and the exact solution. Typically,
this kind of estimates is helpful also in controlling the approximation
error of numerical methods to multiscale problems 
(e.g. \cite{Legoll2014,Hou1999}). The main result of this note is Theorem \ref{thm:1} where we report the upper bound of the corrector up to an arbitrary high order. 

In \cite{KM16}, the investigated microscopic semi-linear system resembles a steady state-type of thermo\textendash diffusion systems (\cite{KMK15,KAM14}). Essentially, we have analyzed the solvability of the microscopic system in \cite{KM16},
derived the upscaled equations as well as the
corresponding effective coefficients, and proved the high-order corrector
estimates for the differences of concentrations and their gradients
in which the standard energy method has been used. Furthermore,
we also solved a reduced similar problem where the Picard iterations-based method
is applied to deal with the nonlinear auxiliary problems (\cite{KM16-1}). 

It is worth noting that the following errors have been so far obtained for the homogenization of the above-mentioned microscopic elliptic system in \cite{KM16} and further in \cite{CP99}:
\[
u^{\varepsilon}-\sum_{m=0}^{M}\varepsilon^{m}u_{m}\quad\text{and}\quad u^{\varepsilon}-u_{0}-m^{\varepsilon}\sum_{m=1}^{M}\varepsilon^{m}u_{m},\quad
M\ge 2.
\]

In this note, we prove a corrector in the form of
\begin{equation}
u^{\varepsilon}-\sum_{k=0}^{K}\varepsilon^{k}u_{k}-m^{\varepsilon}\sum_{m=K+1}^{M}\varepsilon^{m}u_{m},\label{eq:ddd}
\end{equation}
in which we fix $K\in\mathbb{N}$ such that $0\le K\le M-2$.

\section{Problem settings}

We consider the semi-linear elliptic boundary value
problem
\[
\mathcal{A}^{\varepsilon}u_{i}^{\varepsilon}\equiv\nabla\cdot\left(-d_{i}^{\varepsilon}\nabla u_{i}^{\varepsilon}\right)=R_{i}\left(u_{1}^{\varepsilon},...,u_{N}^{\varepsilon}\right)\quad\text{in}\;\Omega^{\varepsilon},
\]
associated with the boundary conditions
\[
d_{i}^{\varepsilon}\nabla u_{i}^{\varepsilon}\cdot\text{n}=\varepsilon\left(a_{i}^{\varepsilon}u_{i}^{\varepsilon}-b_{i}^{\varepsilon}F_{i}\left(u_{i}^{\varepsilon}\right)\right)\quad\text{across}\;\Gamma^{\varepsilon},
\]
\[
u_{i}^{\varepsilon}=0\quad\text{across}\;\Gamma^{ext},
\]
for $i\in\left\{ 1,...,N\right\} $ with $N\ge2$ being the number
of concentrations. For simplicity, we refer this problem as $\left(P^{\varepsilon}\right)$.

\begin{figure}[!h]
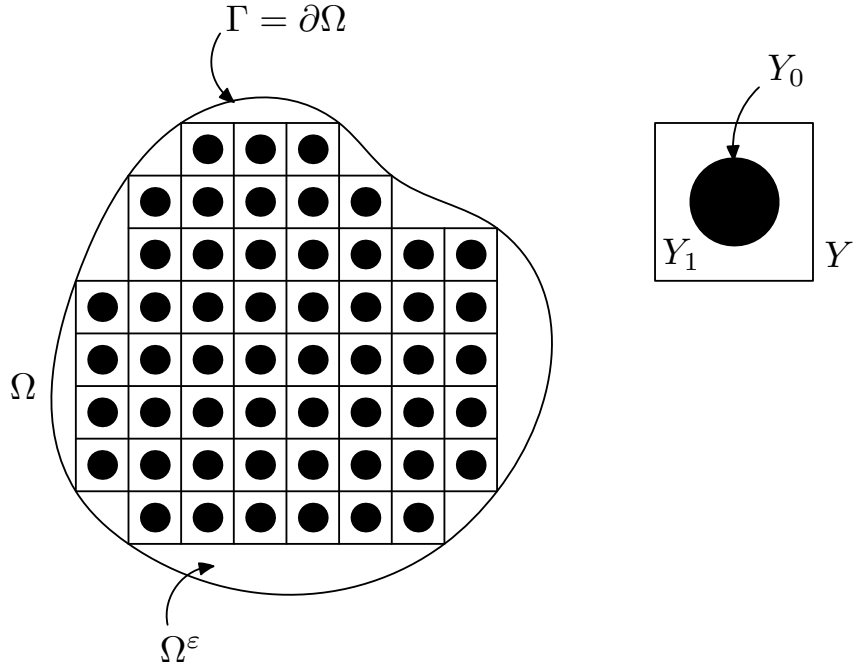
 	
\centering 	
	\parbox{10cm}{\convertMPtoPDF{fig.man}{1.4}{1.4}}	
	\caption{Admissible 2-D perforated domain (left) and basic geometry of the microstructure (right). (By courtesy of Mai Thanh Nhat Truong, Hankyong National University, Republic of Korea.)}
	\label{fig:1}
\end{figure}

This problem is connected to the Smoluchowski--Soret--Dufour modeling of the evolution of temperature and colloid concentrations \cite{dGM62,KMK15}. Here,  $u^{\varepsilon}:=\left(u_{1}^{\varepsilon},...,u_{N}^{\varepsilon}\right)$
denotes the vector of the concentrations, $d_{i}^{\varepsilon}$
represents the molecular diffusion with $R_{i}$ being the volume
reaction rate and $a_{i}^{\varepsilon}$, $b_{i}^{\varepsilon}$ are
deposition coefficients, whilst $F_{i}$ indicates a surface chemical
reaction for the immobile species. Notice that the quantity $\varepsilon$
is called the homogenization parameter or the scale factor. The perforated domain $\Omega^{\varepsilon}\subset\mathbb{R}^{d}$
herein approximates a porous medium and its precise description can be found
in \cite{HJ91,KM16}. As an example, we depict in Figure \ref{fig:1} 
an admissible geometry of our medium and the corresponding microstructure.

Denote by $x\in\Omega^{\varepsilon}$ the macroscopic variable and
by $y=x/\varepsilon$ the microscopic variable representing high oscillations
at the microscopic geometry. Henceforward, we understand throughout
this paper the following convention:
\[
d_{i}^{\varepsilon}\left(x\right)=d_{i}\left(\frac{x}{\varepsilon}\right)=d_{i}\left(y\right),\quad x\in \Omega^{\varepsilon}, y \in Y_{1},
\]
with the same meaning for all the oscillating data such as $a_{i}^{\varepsilon}$,
$b_{i}^{\varepsilon}$, e.g.

We introduce the function space
\[
V^{\varepsilon}:=\left\{ v\in H^{1}\left(\Omega^{\varepsilon}\right)|v=0\;\text{on}\;\Gamma^{ext}\right\} ,
\]
which is a closed subspace of the Hilbert space $H^{1}\left(\Omega^{\varepsilon}\right)$ with the semi-norm
\[
\left\Vert v\right\Vert _{V^{\varepsilon}}=\left(\sum_{i=1}^{d}\int_{\Omega^{\varepsilon}}\left|\frac{\partial v}{\partial x_{i}}\right|^{2}dx\right)^{1/2}\quad\text{for all}\;v\in V^{\varepsilon}.
\]

Denote also $\mathcal{V}^{\varepsilon}=V^{\varepsilon}\times...\times V^{\varepsilon}$ and  $\mathcal{W}^{q}\left(\Omega^{\varepsilon}\right)=L^{q}\left(\Omega^{\varepsilon}\right)\times...\times L^{q}\left(\Omega^{\varepsilon}\right)$ for $q \in (2,\infty]$.

Unless otherwise specified, all the constants $C$ are independent
of the homogenization parameter $\varepsilon$, but the respective values may
differ from line to line and may change even within a single chain
of estimates.

\section{Corrector estimate}

Consider the two--scale asymptotic expansion up to $M$th-level $\left(M\ge2\right)$
given by
\begin{equation}
u_{i}^{\varepsilon}\left(x\right)=\sum_{m=0}^{M}\varepsilon^{m}u_{i,m}\left(x,y\right)+\mathcal{O}\left(\varepsilon^{M+1}\right),\quad x\in\Omega^{\varepsilon},\label{eq:expansion}
\end{equation}
where $u_{i,m}\left(x,\cdot\right)$ is $Y$-periodic for $0\le m\le M$
and $i\in\left\{ 1,...,N\right\} $.

Here we keep the assumptions from \cite{KM16} on
the coefficients. This implicitly guarantees the well-posedness of $\left(P^{\varepsilon}\right)$ 
and $L^{\infty}$ bounds for all concentrations $u_{i}^{\varepsilon}$ for $i\in\left\{1,...,N\right\}$. We also remind the crucial assumptions on the Smoluchowski production term, i.e.
\begin{equation}
R_{i}\left(\sum_{m=0}^{M}\varepsilon^{m}u_{1,m},...,\sum_{m=0}^{M}\varepsilon^{m}u_{N,m}\right)=\sum_{m=0}^{M}\varepsilon^{m}\bar{R}_{i}\left(u_{1,m},...,u_{N,m}\right)+\mathcal{O}\left(\varepsilon^{M+1}\right),\label{assRi}
\end{equation}
\begin{equation}
F_{i}\left(\sum_{m=0}^{M}\varepsilon^{m}u_{i,m}\right)=\sum_{m=0}^{M}\varepsilon^{m}\bar{F}_{i}\left(u_{i,m}\right)+\mathcal{O}\left(\varepsilon^{M+1}\right),\label{assFi}
\end{equation}
where $\bar{R}_{i}$ and $\bar{F}_{i}$ are global Lipschitz functions with the Lipschitz constants $L_{i}$ and $K_{i}$, respectively, for
$i\in\left\{ 1,...,N\right\} $.

Let $m^{\varepsilon}\in C_{c}^{\infty}\left(\Omega\right)$
be a cut-off function such that $\varepsilon\left|\nabla m^{\varepsilon}\right|\le C$
and
\[
m^{\varepsilon}\left(x\right):=\begin{cases}
1, & \text{if}\;\text{dist}\left(x,\Gamma\right)\le\varepsilon,\\
0, & \text{if}\;\text{dist}\left(x,\Gamma\right)\ge2\varepsilon,
\end{cases}
\]
(see \cite{CP99} for more properties of $m_{\varepsilon}$). We aim at proving the following result:

\begin{thm}\label{thm:1}
	Let $u^{\varepsilon}$ be the vector of solutions of the elliptic
	system $\left(P^{\varepsilon}\right)$. Consider the high-order asymptotic
	expansion \eqref{eq:expansion} up to $M$-level $\left(M\ge2\right)$
	and take $\left(u_{0},u_{m}\right)\in\mathcal{W}^{\infty}\left(\Omega^{\varepsilon}\right)\times\mathcal{W}^{\infty}\left(\Omega^{\varepsilon};H_{\#}^{1}\left(Y_{1}\right)/\mathbb{R}\right)$
	for all $0\le m\le M$. For a fixed $K\in\mathbb{N}$ such that $0\le K\le M-2$
	the following corrector estimate holds:
	\begin{equation}
	\left\Vert u^{\varepsilon}-\sum_{k=0}^{K}\varepsilon^{k}u_{k}-m^{\varepsilon}\sum_{m=K+1}^{M}\varepsilon^{m}u_{m}\right\Vert _{\mathcal{V}^{\varepsilon}}\le C\left(\varepsilon^{M-1}+\varepsilon^{M}+\sum_{m=K+1}^{M}\left(\varepsilon^{m-\frac{1}{2}}+\varepsilon^{m+\frac{1}{2}}\right)\right),\label{eq:desired-corrector}
	\end{equation}
	where $C>0$ is a generic $\varepsilon$-independent constant.
\end{thm}

\section{Proof of Theorem \ref{thm:1}}

The following useful
estimates (cf. \cite{Eck2}) hold true:
\begin{equation}
\left\Vert 1-m^{\varepsilon}\right\Vert _{L^{2}\left(\Omega^{\varepsilon}\right)}\le C\varepsilon^{1/2},\quad\varepsilon\left\Vert \nabla m^{\varepsilon}\right\Vert _{L^{2}\left(\Omega^{\varepsilon}\right)}\le C\varepsilon^{1/2}.\label{eq:est-mep}
\end{equation}

To bound from above in term of $\varepsilon$ the quantity \eqref{eq:ddd},
we define the function $\Psi_{i}^{\varepsilon}$ by
\[
\Psi_{i}^{\varepsilon}:=\varphi_{i}^{\varepsilon}+\left(1-m^{\varepsilon}\right)\sum_{m=K+1}^{M}\varepsilon^{m}u_{i,m},
\]
where we denote
\[
\varphi_{i}^{\varepsilon}:=u_{i}^{\varepsilon}-\sum_{m=0}^{M}\varepsilon^{m}u_{i,m}\quad\text{for}\;i\in\left\{ 1,...,N\right\} .
\]

By induction, one can easily obtain that the function $\varphi_{i}^{\varepsilon}$ satisfies
the following equation:
\begin{equation}
\mathcal{A}^{\varepsilon}\varphi_{i}^{\varepsilon}=R_{i}\left(u^{\varepsilon}\right)-\sum_{m=0}^{M-2}\varepsilon^{m}\bar{R}_{i}\left(u_{m}\right)-\varepsilon^{M-1}\left(\mathcal{A}_{1}u_{i,M}+\mathcal{A}_{2}u_{i,M-1}\right)-\varepsilon^{M}\mathcal{A}_{2}u_{i,M}\;\text{in}\;\Omega^{\varepsilon},\label{3.3}
\end{equation}
associated with the following boundary condition at $\Gamma^{\varepsilon}$
\begin{equation}
-d_{i}^{\varepsilon}\nabla_{x}\varphi_{i}^{\varepsilon}\cdot\mbox{n}=\varepsilon^{M}d_{i}^{\varepsilon}\nabla_{x}u_{i,M}\cdot\mbox{n}+\varepsilon\left[a_{i}^{\varepsilon}\left(\sum_{m=0}^{M-2}\varepsilon^{m}u_{i,m}-u_{i}^{\varepsilon}\right)+b_{i}^{\varepsilon}\left(F_{i}\left(u_{i}^{\varepsilon}\right)-\sum_{m=0}^{M-2}\varepsilon^{m}\bar{F}_{i}\left(u_{i,m}\right)\right)\right].\label{3.4}
\end{equation}

In \eqref{3.4}, $\mathcal{A}_{1}$ and $\mathcal{A}_{2}$ are defined, respectively, as follows:
\begin{align*}
\mathcal{A}_{1} & :=\nabla_{x}\cdot\left(-d_{i}\left(y\right)\nabla_{y}\right)+\nabla_{y}\cdot\left(-d_{i}\left(y\right)\nabla_{x}\right),\\
\mathcal{A}_{2} & :=\nabla_{x}\cdot\left(-d_{i}\left(y\right)\nabla_{x}\right).
\end{align*}

Multiplying \eqref{3.3} by $\varphi_{i}\in V^{\varepsilon}$ and
integrating by parts with using \eqref{3.4}, we arrive at
\begin{eqnarray}
\int_{\Omega^{\varepsilon}}d_{i}^{\varepsilon}\nabla\varphi_{i}^{\varepsilon}\nabla\varphi_{i}dx & = & \left\langle R_{i}\left(u^{\varepsilon}\right)-\sum_{m=0}^{M-2}\varepsilon^{m}\bar{R}_{i}\left(u_{m}\right),\varphi_{i}\right\rangle _{L^{2}\left(\Omega^{\varepsilon}\right)}\nonumber \\
 &  & -\varepsilon^{M-1}\left\langle \mathcal{A}_{1}u_{i,M}+\mathcal{A}_{2}u_{i,M-1}+\varepsilon\mathcal{A}_{2}u_{i,M},\varphi_{i}\right\rangle _{L^{2}\left(\Omega^{\varepsilon}\right)}\nonumber \\
 &  & -\varepsilon\left\langle a_{i}^{\varepsilon}\left(\sum_{m=0}^{M-2}\varepsilon^{m}u_{i,m}-u_{i}^{\varepsilon}\right)+b_{i}^{\varepsilon}\left(F_{i}\left(u_{i}^{\varepsilon}\right)-\sum_{m=0}^{M-2}\varepsilon^{m}\bar{F}_{i}\left(u_{i,m}\right)\right),\varphi_{i}\right\rangle _{L^{2}\left(\Gamma^{\varepsilon}\right)}\nonumber \\
 &  & -\varepsilon^{M}\int_{\Gamma^{\varepsilon}}d_{i}^{\varepsilon}\nabla_{x}u_{i,M}\cdot\mbox{n}\varphi_{i}dS_{\varepsilon}.\label{eq:4.29}
\end{eqnarray}

We can now gain the first part of the corrector
\eqref{eq:desired-corrector}, 
i.e. we shall estimate each integral on the right-hand side of \eqref{eq:4.29} which we denote by $\mathcal{I}_1$, $\mathcal{I}_2$, $\mathcal{I}_3$ and $\mathcal{I}_4$, respectively.

Let $\bar{L}:=\max\left\{ \bar{L}_{1},...,\bar{L}_{N}\right\} $.
Using \eqref{assRi} in combination with the structural inequality $\left\Vert \bar{R}_{i}\left(u_{m}\right)\right\Vert _{L^{2}\left(\Omega^{\varepsilon}\right)}\le\bar{L}\left\Vert u_{m}\right\Vert _{\mathcal{W}^{2}\left(\Omega^{\varepsilon}\right)}+\left\Vert \bar{R}_{i}\left(0\right)\right\Vert _{L^{2}\left(\Omega^{\varepsilon}\right)}$
for all $0\le m\le M$, we see that
\begin{align}
\left|\left\langle R_{i}\left(u^{\varepsilon}\right)-\sum_{m=0}^{M-2}\varepsilon^{m}\bar{R}_{i}\left(u_{m}\right),\varphi_{i}\right\rangle _{L^{2}\left(\Omega^{\varepsilon}\right)}\right| & \le\varepsilon^{M-1}\left(\bar{L}\left\Vert u_{M-1}\right\Vert _{\mathcal{W}^{2}\left(\Omega^{\varepsilon}\right)}+\left\Vert \bar{R}_{i}\left(0\right)\right\Vert _{L^{2}\left(\Omega^{\varepsilon}\right)}\right)\left\Vert \varphi_{i}\right\Vert _{L^{2}\left(\Omega^{\varepsilon}\right)}\nonumber \\
 & +\varepsilon^{M}\left(\bar{L}\left\Vert u_{M}\right\Vert _{\mathcal{W}^{2}\left(\Omega^{\varepsilon}\right)}+\left\Vert \bar{R}_{i}\left(0\right)\right\Vert _{L^{2}\left(\Omega^{\varepsilon}\right)}\right)\left\Vert \varphi_{i}\right\Vert _{L^{2}\left(\Omega^{\varepsilon}\right)}\nonumber \\
 & \le C\left(\varepsilon^{M-1}+\varepsilon^{M}\right)\left\Vert \varphi_{i}\right\Vert _{L^{2}\left(\Omega^{\varepsilon}\right)}.\label{eq:2.15}
\end{align}

The second integral $\mathcal{I}_2$ can be bounded above by
\begin{equation}
\varepsilon^{M-1}\left|\left\langle \mathcal{A}_{1}u_{i,M}+\mathcal{A}_{2}u_{i,M-1}+\varepsilon\mathcal{A}_{2}u_{i,M},\varphi_{i}\right\rangle _{L^{2}\left(\Omega^{\varepsilon}\right)}\right|\le C\varepsilon^{M-1}\left\Vert \varphi_{i}\right\Vert _{L^{2}\left(\Omega^{\varepsilon}\right)}.\label{eq:2.16}
\end{equation}

Let $\bar{K}:=1+\max\left\{ \bar{K}_{1},...,\bar{K}_{N}\right\} $.
For the integral $\mathcal{I}_3$, we proceed as the proof of 
\eqref{eq:2.15}. We thus claim that
\begin{equation}
\varepsilon\left|\left\langle a_{i}^{\varepsilon}\left(\sum_{m=0}^{M-2}\varepsilon^{m}u_{i,m}-u_{i}^{\varepsilon}\right)+b_{i}^{\varepsilon}\left(F_{i}\left(u_{i}^{\varepsilon}\right)-\sum_{m=0}^{M-2}\varepsilon^{m}\bar{F}_{i}\left(u_{i,m}\right)\right),\varphi_{i}\right\rangle _{L^{2}\left(\Gamma^{\varepsilon}\right)}\right|\le C\left(\varepsilon^{M-1}+\varepsilon^{M}\right)\left\Vert \varphi_{i}\right\Vert _{L^{2}\left(\Omega^{\varepsilon}\right)},\label{3.10}
\end{equation}
in which we have used \eqref{eq:expansion} and \eqref{assFi} together
with the H\"older inequality, as well as the trace inequality (cf. \cite[Lemma 2.31]{CP99}). On top of that, it yields for the last integral $\mathcal{I}_4$ that
\begin{align}
\varepsilon^{M}\left|\int_{\Gamma^{\varepsilon}}d_{i}^{\varepsilon}\nabla_{x}u_{i,M}\cdot\mbox{n}\varphi_{i}dS_{\varepsilon}\right| & \le\varepsilon^{M}\left\Vert d_{i}^{\varepsilon}\nabla_{x}u_{i,M}\cdot\mbox{n}\right\Vert _{L^{2}\left(\Gamma^{\varepsilon}\right)}\left\Vert \varphi_{i}\right\Vert _{L^{2}\left(\Gamma^{\varepsilon}\right)}\nonumber \\
 & \le C\varepsilon^{M-1}\left\Vert \varphi_{i}\right\Vert _{L^{2}\left(\Omega^{\varepsilon}\right)}.\label{eq:2.17}
\end{align}
Combining \eqref{eq:4.29}-\eqref{eq:2.17}, we observe
\begin{equation}
\left|\left\langle \varphi_{i}^{\varepsilon},\varphi_{i}\right\rangle _{V^{\varepsilon}}\right|\le C\left(\varepsilon^{M-1}+\varepsilon^{M}\right)\left\Vert \varphi_{i}\right\Vert _{L^{2}\left(\Omega^{\varepsilon}\right)}\quad\text{for}\;\varphi_{i}\in V^{\varepsilon}\;\text{and}\;i\in\left\{ 1,...,N\right\} ,\label{err1}
\end{equation}
which then leads to $\left\Vert \varphi_{i}^{\varepsilon}\right\Vert _{V^{\varepsilon}}\le C\varepsilon^{M-1}$
by choosing $\varphi_{i}=\varphi_{i}^{\varepsilon}$ for $i\in\left\{ 1,...,N\right\} $.

It remains to consider the following quantity:
\[
\left\langle \left(1-m^{\varepsilon}\right)\sum_{m=K+1}^{M}\varepsilon^{m}u_{i,m},\varphi_{i}\right\rangle _{V^{\varepsilon}}\quad\text{for}\;\varphi_{i}\in V^{\varepsilon}\;\text{and}\;i\in\left\{ 1,...,N\right\} .
\]

At this stage, the following estimate is straightforward due to \eqref{eq:est-mep}:
\begin{align}
\left|\left\langle \left(1-m^{\varepsilon}\right)\sum_{m=K+1}^{M}\varepsilon^{m}u_{i,m},\varphi_{i}\right\rangle _{V^{\varepsilon}}\right| & \le C\left\Vert \nabla\left(1-m^{\varepsilon}\right)\left(\sum_{m=K+1}^{M}\varepsilon^{m}u_{i,m}\right)\right\Vert _{L^{2}\left(\Omega^{\varepsilon}\right)}\left\Vert \varphi_{i}\right\Vert _{V^{\varepsilon}}\nonumber \\
 & +C\left\Vert \left(1-m^{\varepsilon}\right)\nabla\left(\sum_{m=K+1}^{M}\varepsilon^{m}u_{i,m}\right)\right\Vert _{L^{2}\left(\Omega^{\varepsilon}\right)}\left\Vert \varphi_{i}\right\Vert _{V^{\varepsilon}}\nonumber \\
 & \le C\sum_{m=K+1}^{M}\varepsilon^{m}\left\Vert \nabla\left(1-m^{\varepsilon}\right)\right\Vert _{L^{2}\left(\Omega^{\varepsilon}\right)}\left\Vert \varphi_{i}\right\Vert _{V^{\varepsilon}}\nonumber \\
 & +C\sum_{m=K+1}^{M}\varepsilon^{m}\left\Vert 1-m^{\varepsilon}\right\Vert _{L^{2}\left(\Omega^{\varepsilon}\right)}\left\Vert \varphi_{i}\right\Vert _{V^{\varepsilon}}\nonumber \\
 & \le C\sum_{m=K+1}^{M}\left(\varepsilon^{m-\frac{1}{2}}+\varepsilon^{m+\frac{1}{2}}\right)\left\Vert \varphi_{i}\right\Vert _{V^{\varepsilon}}\; \mbox{for all}\;\varphi_{i}\in V^{\varepsilon}.\label{eq:err2}
\end{align}

Thanks to the triangle inequality, we combine \eqref{err1} and \eqref{eq:err2}
and then sum up the resulting estimates to
get
\[
\sum_{i=1}^{N}\left|\left\langle \Psi_{i}^{\varepsilon},\varphi_{i}\right\rangle _{V^{\varepsilon}}\right|\le C\left(\varepsilon^{M-1}+\varepsilon^{M}+\sum_{m=K+1}^{M}\left(\varepsilon^{m-\frac{1}{2}}+\varepsilon^{m+\frac{1}{2}}\right)\right)\left\Vert \varphi\right\Vert _{\mathcal{V}^{\varepsilon}}\quad\text{for}\;\varphi\in\mathcal{V}^{\varepsilon}.
\]

By choosing $\varphi=\Psi^{\varepsilon}$ and then by simplifying both
sides of the resulting estimate by $\left\Vert \Psi^{\varepsilon}\right\Vert _{\mathcal{V}^{\varepsilon}}$,
we obtain that
\[
\left\Vert \Psi^{\varepsilon}\right\Vert _{\mathcal{V}^{\varepsilon}}\le C\left(\varepsilon^{M-1}+\varepsilon^{M}+\sum_{m=K+1}^{M}\left(\varepsilon^{m-\frac{1}{2}}+\varepsilon^{m+\frac{1}{2}}\right)\right).
\]

This completes the proof of Theorem \ref{thm:1}.

\section*{Acknowledgements}
Helpful discussions with Prof. Adrian Muntean are acknowledged.

\bibliographystyle{plain}
\bibliography{mybib}

\begin{thebibliography}{10}

\bibitem{Legoll2014}
C.~Le Bris, F.~Legoll, and A.~Lozinski.
\newblock An {M}s{FEM} type approach for perforated domains.
\newblock {\em SIAM Multiscale Modeling and Simulation}, 12(3):1046--1077,
  2014.

\bibitem{CP99}
D.~Cioranescu and {J. Saint Jean Paulin}.
\newblock {\em Homogenization of {R}eticulated {S}tructures}.
\newblock Springer, 1999.

\bibitem{dGM62}
S.R. de~Groot and P.~Mazur.
\newblock {\em Non-equilibrium Thermodynamics}.
\newblock North-Holland Publishing Company, Amsterdam, 1962.

\bibitem{Eck2}
C.~Eck.
\newblock Homogenization of a phase field model for binary mixtures.
\newblock {\em Multiscale Modeling and Simulation}, 3:1--27, 2004.

\bibitem{HJ91}
U.~Hornung and W.~J\"ager.
\newblock {D}iffusion, convection, adsorption, and reaction of chemicals in
  porous media.
\newblock {\em Journal of Differential Equations}, 92:199--225, 1991.

\bibitem{Hou1999}
T.~Hou, X.~H. Wu, and Z.~Cai.
\newblock Convergence of a multiscale finite element method for elliptic
  problems with rapidly oscillating coefficients.
\newblock {\em Mathematics of Computation}, 68(227):913--943, 1999.

\bibitem{KM16}
V.A. Khoa and A.~Muntean.
\newblock Asymptotic analysis of a semi-linear elliptic system in perforated
  domains: Well-posedness and corrector for the homogenization limit.
\newblock {\em Journal of Mathematical Analysis and Applications},
  439:271--295, 2016.

\bibitem{KM16-1}
V.A. Khoa and A.~Muntean.
\newblock A note on iterations-based derivations of high-order homogenization
  correctors for multiscale semi-linear elliptic equations.
\newblock {\em Applied Mathematics Letters}, 58:103--109, 2016.

\bibitem{KAM14}
O.~Krehel, T.~Aiki, and A.~Muntean.
\newblock Homogenization of a thermo-diffusion system with {S}moluchowski
  interactions.
\newblock {\em Networks and Heterogeneous Media}, 9(4):739--762, 2014.

\bibitem{KMK15}
O.~Krehel, A.~Muntean, and P.~Knabner.
\newblock Multiscale modeling of colloidal dynamics in porous media including
  aggregation and deposition.
\newblock {\em Advances in Water Resources}, 86:209--216, 2015.

\end{thebibliography}

\end{document}